\documentclass{amsart}
\usepackage{amsmath}
\usepackage{amssymb}
\usepackage{amsfonts}

\theoremstyle{plain}

\numberwithin{equation}{section}
\input{tcilatex}

\begin{document}
\title[On Orthogonal Decomposition of a Sobolev Space]{On Orthogonal
Decomposition of a Sobolev Space}
\author{Dejenie A. Lakew}
\address{Bryant \& Stratton College\\
Richmond, VA 23235}
\email{dalakew@bryantstratton.edu}
\urladdr{www.bryantstratton.edu}
\date{November 13, 2016}

\begin{abstract}
The theme of this short article is to investigate an orthogonal
decomposition of the Sobolev space $W^{1,2}\left( \Omega \right) $ as $%
W^{1,2}\left( \Omega \right) =A^{2,2}\left( \Omega \right) \oplus
D^{2}\left( W_{0}^{3,2}\left( \Omega \right) \right) $ and look at some
properties of the inner product therein and the distance defined from the
inner product. We also determine the dimension of\ the orthogonal difference
space $W^{1,2}\left( \Omega \right) \ominus \left( W_{0}^{1,2}\left( \Omega
\right) \right) ^{\perp }$ and show the expansion of Sobolev spaces as their
regularity increases.
\end{abstract}

\maketitle

\section{\protect\bigskip}

This is an extension work of $[1]$ and $\left[ 2\right] $ in which the space
under consideration is a Sobolev space of regularity exponent one. The
change in regularity, from the Lebesgue space of regularity zero to Sobolev
spaces of higher regularities, causes increase in length or norm, expansion
of the space in terms of distance or separation between distinct elements
and change in orthogonality.

\ \ \ \ \ \ 

In addition to the regular properties we develop in the decomposition
process, we obtain some geometric properties of distance and apertures as
well between non zero elements.

\ \ 

\textbf{Notations.}

\ 

$\Omega :=[0,1]$

\ 

$\oplus :=$ Direct sum for sets, $\ominus :=$ Direct difference of sets

\ 

$\uplus :=$ Direct sum of elements from orthogonal sets

\ 

$D^{\alpha }:=\frac{d^{\alpha }}{dx^{\alpha }}$, for $\alpha =0,1,2$

$C^{\infty }\left( \Omega \right) =\dbigcap\limits_{n=0}^{\infty
}C^{n}\left( \Omega \right) $

$C_{0}^{\infty }\left( \Omega \right) =\left\{ f\in C^{\infty }\left( \Omega
\right) :\sup pf\subset K\Subset \Omega \right\} $\ \ \ 

\ \ \ \ 

\textbf{Definition 1}. We say that a function $h:\Omega \longrightarrow 
\mathbb{R}
$ is a weak derivative of $g$ of order $\alpha $ if

\ 

\begin{equation*}
\dint\limits_{\Omega }h(x)\phi (x)dx=\left( -1\right) ^{\alpha
}\dint\limits_{\Omega }g(x)D^{\alpha }\phi \left( x\right) dx,\forall \phi
\in C_{0}^{\infty }\left( \Omega \right) .
\end{equation*}

\ 

Clearly functions which are differentiable in the regular sense of a certain
order are weakly differentiable of that order but not the converse.

\ \ \ 

\textbf{Example 1}.%
\begin{equation*}
f(x)=\left\{ 
\begin{array}{c}
x-\frac{1}{2},\frac{1}{2}\leq x\leq 1 \\ 
0,0\leq x\leq \frac{1}{2}%
\end{array}%
\right.
\end{equation*}

Then%
\begin{equation*}
f\in C^{0}\left( \Omega \right) \backslash C^{1}(\Omega )
\end{equation*}

\ \ \ \ \ \ \ 

i.e. $f$ is continuous but not differentiable in the regular sense but
weakly differentiable with weak derivative

\ \ \ 

\begin{equation*}
Df=g=\left\{ 
\begin{array}{c}
1,\frac{1}{2}<x\leq 1 \\ 
0,0\leq x<\frac{1}{2}%
\end{array}%
\right.
\end{equation*}

Indeed $\ $%
\begin{eqnarray*}
\dint\limits_{\Omega }f\phi ^{\prime }dx &=&\dint\limits_{\Omega }\left( x-%
\frac{1}{2}\right) \phi ^{\prime }dx \\
&=&-\dint \left( x-\frac{1}{2}\right) ^{\prime }\phi dx \\
&=&-\dint\limits_{\frac{1}{2}}^{1}\phi dx \\
&=&-\dint\limits_{0}^{1}Df\phi dx
\end{eqnarray*}

\textbf{Definition 2}. The Sobolev space 
\begin{equation*}
W^{1,2}\left( \Omega \right) :=\{f\in \tciLaplace ^{2}\left( \Omega \right)
:Df\in \tciLaplace ^{2}(\Omega )\}
\end{equation*}

and \ \ \ \ 
\begin{equation*}
W_{0}^{1,2}\left( \Omega \right) =\left\{ f\in W^{1,2}\left( \Omega \right)
:f_{\mid \partial \Omega }=0\right\}
\end{equation*}

where $Df$ \ is in the sense of weak/distributional derivative.

\ \ \ \ \ \ \ \ \ \ \ \ \ \ \ \ \ \ \ \ \ \ \ \ \ \ \ \ 

\textbf{Note}. 
\begin{equation*}
W^{1,2}\left( \Omega \right) \subseteq \tciLaplace ^{2}\left( \Omega \right)
\end{equation*}

but not the converse, 
\begin{equation*}
\tciLaplace ^{2}\left( \Omega \right) \nsubseteq W^{1,2}\left( \Omega \right)
\end{equation*}

\ 

\textbf{Example 2}. The function%
\begin{equation*}
f\left( x\right) =\sqrt{x}\in \tciLaplace ^{2}\left( \Omega \right)
\end{equation*}

but 
\begin{equation*}
f\notin W^{1,2}\left( \Omega \right)
\end{equation*}

\ 

since $f^{\prime }\left( x\right) $ has a singularity at $x=0$ where the
improper integral

\begin{equation*}
\dint\limits_{\Omega }\mid f^{\prime }\left( x\right) \mid ^{2}dx
\end{equation*}

diverges. That is $f^{\prime }\notin \tciLaplace ^{2}\left( \Omega \right) $

\ \ 

\textbf{Question}: Is there a non trivial function%
\begin{equation*}
f\in W^{1,2}\left( \Omega \right) \ominus W_{0}^{1,2}\left( \Omega \right)
\end{equation*}

and how big is 
\begin{equation*}
W^{1,2}\left( \Omega \right) \ominus W_{0}^{1,2}\left( \Omega \right) \text{
\ \ ?}
\end{equation*}

\ 

\textbf{Embedding}. Clearly the Sobolev space $W^{1,2}\left( \Omega \right) $
is not a collection of wildly behaved generalized functions but some how
well behaved functions that are more than continuous.

\ 

In fact the space is embedded in the H\"{o}lder space $C^{\gamma }$ for $%
0\leq \gamma \leq \frac{1}{2}$, in particularly $W^{1,2}\left( \Omega
\right) \hookrightarrow C^{\frac{1}{2}}\left( \Omega \right) $. The H\"{o}%
lder space $C^{\frac{1}{2}}\left( \Omega \right) $ is the space of functions 
$f$ with property

\ \ \ 
\begin{equation*}
\parallel f\left( x\right) -f\left( y\right) \parallel \leq \lambda
_{f}\parallel x-y\parallel ^{\frac{1}{2}},\forall x,y\in \Omega
\end{equation*}

\ for some non negative constant $\lambda _{f}$ that depends on $f$.

\ \ \ \ \ \ \ \ 

\textbf{Inner Product}. The Sobolev space $W^{1,2}\left( \Omega \right) $ is
an inner product space with inner product defined as

\ 

\begin{equation*}
\left\langle f,g\right\rangle _{W^{1,2}\left( \Omega \right)
}:=\dint\limits_{\Omega }fg+f^{\prime }g^{\prime }dx
\end{equation*}

and from this inner product we define norm as

\begin{equation*}
\parallel f\parallel _{W^{1,2}\left( \Omega \right) }=\left( \left\langle
f,f\right\rangle _{W^{1,2}\left( \Omega \right) }\right) ^{\frac{1}{2}}
\end{equation*}

Clearly $\ $

\ \ \ \ \ \ \ \ \ \ 

$\left( i\right) $ 
\begin{equation*}
\parallel \bullet \parallel _{W^{1,2}\left( \Omega \right) }\geq \parallel
\bullet \parallel _{\tciLaplace ^{2}\left( \Omega \right) }
\end{equation*}

\ \ \ \ \ \ \ \ \ \ \ \ \ 

$\left( ii\right) \ \ $%
\begin{equation*}
\left\langle f,g\right\rangle _{W^{1,2}\left( \Omega \right) }=\left\langle
f,g\right\rangle _{\tciLaplace ^{2}\left( \Omega \right) }+\left\langle
f^{\prime },g^{\prime }\right\rangle _{\tciLaplace ^{2}\left( \Omega \right)
}
\end{equation*}

\ \ \ 

\textbf{Definition 3}. Two functions $f$ and $g$ of $W^{1,2}\left( \Omega
\right) $ are said to be orthogonal with respect to the inner product
defined above iff

\ \ 
\begin{equation*}
\left\langle f,g\right\rangle _{W^{1,2}\left( \Omega \right) }=0
\end{equation*}

\ \ \ 

\textbf{Example 3}.

\ 

$\left( a\right) $ 
\begin{equation*}
\left\langle \sin x,\sin x\right\rangle _{_{W^{1,2}\left( \Omega \right)
}}=\left( \dint\limits_{\Omega }\sin ^{2}x+\cos ^{2}xdx\right) ^{\frac{1}{2}%
}=1
\end{equation*}

Hence 
\begin{equation*}
\parallel \sin x\parallel _{W^{1,2}\left( \Omega \right) }=1
\end{equation*}

\ 

$\left( b\right) $ 
\begin{equation*}
\left\langle \sin x,\cos x\right\rangle _{_{W^{1,2}\left( \Omega \right) }}=0
\end{equation*}

\ \ \ 

and hence $\sin x$ and $\cos x$ are orthogonal in $W^{1,2}\left( \Omega
\right) $ but not in $\tciLaplace ^{2}\left( \Omega \right) $

\ \ 

$\left( c\right) $ 
\begin{equation*}
\left\langle e^{\alpha x},e^{\beta x}\right\rangle _{_{W^{1,2}\left( \Omega
\right) }}=0\text{ \ for }\alpha \beta =-1
\end{equation*}

\ \ \ \ \ \ \ \ \ \ \ \ \ \ \ \ \ \ \ \ \ \ \ \ \ \ \ \ \ \ \ \ \ \ \ \ 

\ but not in $\left\langle .,.\right\rangle _{\tciLaplace ^{2}\left( \Omega
\right) }$ and hence for $\alpha \beta =-1$ , 
\begin{equation*}
e^{\alpha x}\text{ \ \ and \ \ }e^{\beta x}
\end{equation*}

are orthogonal in $W^{1,2}\left( \Omega \right) $ not in $\tciLaplace
^{2}\left( \Omega \right) $.

\ \ \ 

In particular for $\alpha =1,\beta =-1$, 
\begin{equation*}
\left\langle e^{x},e^{-x}\right\rangle _{_{W^{1,2}\left( \Omega \right) }}=0
\end{equation*}

\ 

\textbf{Proposition 1}. $\ $%
\begin{equation*}
\left\langle f,\lambda f\right\rangle _{_{W^{1,2}\left( \Omega \right)
}}=\lambda \parallel f\parallel _{W^{1,2}\left( \Omega \right) }^{2}
\end{equation*}

\ 

\textbf{Proof}.

\begin{eqnarray*}
\left\langle f,\lambda f\right\rangle _{_{W^{1,2}\left( \Omega \right) }}
&=&\dint\limits_{\Omega }\lambda f^{2}+\lambda \left( f^{\prime }\right)
^{2}dx \\
&=&\lambda \left( \dint\limits_{\Omega }f^{2}+\left( f^{\prime }\right)
^{2}\right) dx \\
&=&\lambda \parallel f\parallel _{W^{1,2}\left( \Omega \right) }^{2}
\end{eqnarray*}

\ \ \ \ \ \ 

\textbf{Proposition 2}. Norm is longer in $W^{1,2}\left( \Omega \right) $
than in $\tciLaplace ^{2}\left( \Omega \right) $, i.e.

\ \ \ \ \ \ \ \ \ \ \ 

$\left( a\right) $ For $f\in W^{1,2}\left( \Omega \right) $ we have

\ \ \ \ \ 
\begin{equation*}
\parallel f\parallel _{W^{1,2}\left( \Omega \right) }\geq \parallel
f\parallel _{\tciLaplace ^{2}\left( \Omega \right) }
\end{equation*}

\ 

$\left( b\right) $ For $f\in W^{1,2}\left( \Omega \right) \cap C^{1}\left(
\Omega \right) :f^{\prime }=\alpha f$ for $\alpha \neq 0$ we have

\ \ \ \ \ \ 
\begin{equation*}
\parallel f\parallel _{W^{1,2}\left( \Omega \right) }=\sqrt{1+\alpha ^{2}}%
\parallel f\parallel _{\tciLaplace ^{2}\left( \Omega \right) }
\end{equation*}

\textbf{Proof}. \ 

\ \ \ 

$\left( a\right) $ \ Let $f\in W^{1,2}\left( \Omega \right) $. Then 
\begin{eqnarray*}
&\parallel &f\parallel _{W^{1,2}\left( \Omega \right) }=\left(
\dint\limits_{\Omega }f^{2}+f^{\prime 2}dx\right) ^{\frac{1}{2}} \\
&\geq &\left( \dint\limits_{\Omega }f^{2}dx\right) ^{\frac{1}{2}} \\
&=&\text{ }\parallel f\parallel _{\tciLaplace ^{2}\left( \Omega \right) }
\end{eqnarray*}

\ \ \ 

$\left( b\right) \ $For $f\in W^{1,2}\left( \Omega \right) \cap C^{1}\left(
\Omega \right) :f^{\prime }=\alpha f$ for $\alpha \neq 0$ , then

\begin{eqnarray*}
&\parallel &f\parallel _{W^{1,2}\left( \Omega \right) }=\left(
\dint\limits_{\Omega }f^{2}+f^{\prime 2}dx\right) ^{\frac{1}{2}} \\
&=&\left( \dint\limits_{\Omega }f^{2}+\left( \alpha f\right) ^{2}dx\right) ^{%
\frac{1}{2}} \\
&=&\left( \dint\limits_{\Omega }f^{2}+\alpha ^{2}f^{2}dx\right) ^{\frac{1}{2}%
} \\
&=&\left( \dint\limits_{\Omega }\left( 1+\alpha ^{2}\right) f^{2}dx\right) ^{%
\frac{1}{2}} \\
&=&\sqrt{1+\alpha ^{2}}\parallel f\parallel _{\tciLaplace ^{2}\left( \Omega
\right) }
\end{eqnarray*}%
\ \ 

Note here that $\left( a\right) $ follows from $\left( b\right) $ since $%
\sqrt{1+\alpha ^{2}}\geq 1$

\ \ \ \ \ \ \ \ \ 

An important question one can pause: which elements $f\in W^{1,2}\left(
\Omega \right) $ maintain their norms of $\tciLaplace ^{2}\left( \Omega
\right) $ ? The answer lies in the next proposition

\ \ \ 

\textbf{Proposition 3.} A function $f\in W^{1,2}\left( \Omega \right) $
which is a.e. a constant over $\Omega $ maintains its norm of $\tciLaplace
^{2}\left( \Omega \right) $.

\ 

\textbf{Proof}. Let $f\in W^{1,2}\left( \Omega \right) $ which is a.e. a
constant. That is $f^{\prime }=0$ a.e. over $\Omega $. Then

\ \ 

\begin{equation*}
\parallel f\parallel _{W^{1,2}\left( \Omega \right) }=\sqrt{\langle
f,f\rangle _{W^{1,2}\left( \Omega \right) }}=\sqrt{\dint\limits_{\Omega
}\left( f^{2}+f^{\prime 2}\right) dx}
\end{equation*}
\ \ and

\begin{equation*}
\sqrt{\dint\limits_{\Omega }\left( f^{2}+f^{\prime 2}\right) dx}=\parallel
f\parallel _{\tciLaplace ^{2}\left( \Omega \right) }
\end{equation*}
only when%
\begin{equation*}
\dint\limits_{\Omega }f^{\prime 2}dx=0
\end{equation*}
That is, $f^{\prime }=0$ a.e. over $\Omega $ means $f$ is a constant a.e.

\ \ \ 

\textbf{Example 4}.

\ \ 

$\left( a\right) $ 
\begin{equation*}
\parallel e^{\pm x}\parallel _{W^{1,2}\left( \Omega \right) }=\sqrt{2}%
\parallel e^{\pm x}\parallel _{\tciLaplace ^{2}\left( \Omega \right) }
\end{equation*}

\ 

$\left( b\right) $ 
\begin{equation*}
\parallel x\parallel _{W^{1,2}\left( \Omega \right) }=2\parallel x\parallel
_{\tciLaplace ^{2}\left( \Omega \right) }
\end{equation*}

\ \ 

In the following propositon, we show that if a function $f\in C^{1}\left(
\Omega \right) $ and its derivative $f^{\prime }$ have a vanishing property
over the boundary of $\Omega $, then always $f$ and its derivative $%
f^{\prime }$ are orthogonal over $\Omega $.

\ \ \ \ \ \ \ \ \ \ \ \ \ \ 

\textbf{Proposition 4}. \ For $f\in C^{1}\left( \Omega \right) $ with
boundary conditions%
\begin{equation*}
f_{\mid \partial \Omega }=0,\text{ \ \ }f_{\mid \partial \Omega }^{\prime }=0
\end{equation*}

then $f$ and $f^{\prime }$ are orthogonal in $W^{1,2}\left( \Omega \right) $%
, i.e. 
\begin{equation*}
\left\langle f,f^{\prime }\right\rangle _{W^{1,2}\left( \Omega \right) }=0
\end{equation*}

\textbf{Proof}. Let $f\in C^{1}\left( \Omega \right) $ than%
\begin{equation*}
\left\langle f,f^{\prime }\right\rangle _{_{W^{1,2}\left( \Omega \right)
}}=\dint\limits_{\Omega }ff^{\prime }+f^{\prime }f^{\prime \prime }dx
\end{equation*}

But%
\begin{eqnarray*}
ff^{\prime }+f^{\prime }f^{\prime \prime } &=&\frac{1}{2}\left( \left(
f^{2}\right) ^{\prime }+\left( \left( f^{\prime }\right) ^{2}\right)
^{\prime }\right) \\
&=&\frac{1}{2}\left( f^{2}+\left( f^{\prime }\right) ^{2}\right) ^{\prime }
\end{eqnarray*}

Thus 
\begin{eqnarray*}
\dint\limits_{\Omega }ff^{\prime }+f^{\prime }f^{\prime \prime }dx &=&\frac{1%
}{2}\dint\limits_{\Omega }\left( f^{2}+\left( f^{\prime }\right) ^{2}\right)
^{\prime }dx \\
&=&\frac{1}{2}\left( f^{2}+\left( f^{\prime }\right) ^{2}\right) \mid
_{0}^{1}\text{ }=0
\end{eqnarray*}

\ \ \ \ \ \ \ \ \ \ \ \ 

\textbf{Remark}. The converse of Proposition 3 does not hold true, since $%
\sin x$ and $\cos x$ are $W^{1,2}\left( \Omega \right) -$ orthogonal but $%
\sin x_{\mid \partial \Omega }\neq 0$ and $\cos x_{\mid \partial \Omega
}\neq 0$.

\ \ 

\ 

\textbf{Example 5}. For $\alpha >1,\beta >1$ let

\ \ \ 
\begin{equation*}
f(x)=x^{\alpha }\left( x-1\right) ^{\beta }
\end{equation*}

then from the above proposition we have

\ 
\begin{equation*}
\left\langle f,f^{\prime }\right\rangle _{W^{1,2}\left( \Omega \right) }=0
\end{equation*}

\ 

\textbf{Proposition 5}.

\ \ \ \ \ 

$\left( a\right) $ 
\begin{equation*}
\left\langle f,g\right\rangle _{W^{1,2}\left( \Omega \right) }=0\text{ \ }%
\Longrightarrow \text{ \ }\left\langle f,g\right\rangle _{\tciLaplace
^{2}\left( \Omega \right) }=-\left\langle f^{\prime },g^{\prime
}\right\rangle _{\tciLaplace ^{2}\left( \Omega \right) }
\end{equation*}

\ \ \ \ 

$\left( b\right) $ 
\begin{equation*}
\left\langle f,g\right\rangle _{\tciLaplace ^{2}\left( \Omega \right) }=0%
\text{ \ }\Longrightarrow \text{ \ }\left\langle f,g\right\rangle
_{W^{1,2}\left( \Omega \right) }=\left\langle f^{\prime },g^{\prime
}\right\rangle _{\tciLaplace ^{2}\left( \Omega \right) }
\end{equation*}

\ \ \ \ 

$(c)$ \ When pair wise $f$ and $g$ and $f^{\prime }$and $g^{\prime }$ are $%
\tciLaplace ^{2}\left( \Omega \right) -$ orthogonal then $f$ and $g$ are $%
W^{1,2}\left( \Omega \right) -$ orthogonal.

\ \ \ 

\ \ \ 

\ \ 

\textbf{Proposition 6. \ }Let $f$ be a function in $W^{1,2}\left( \Omega
\right) $ with non zero norm and $\alpha $ be a non zero constant. Then

\ \ 
\begin{equation*}
\cos \theta \langle f,\alpha f)=\left\{ 
\begin{array}{c}
1\text{, for }\alpha >0 \\ 
-1\text{, for }\alpha <0%
\end{array}%
\right.
\end{equation*}

\ \ \ \ 

\textbf{Proof}. $\ $%
\begin{eqnarray*}
\cos \theta \langle \left( f,\alpha f\right) _{W^{1,2}\left( \Omega \right)
} &=&\frac{\left\langle f,\alpha f\right\rangle _{W^{1,2}\left( \Omega
\right) }}{\mid \alpha \mid \parallel f\parallel _{W^{1,2}\left( \Omega
\right) }^{2}} \\
&=&\frac{\alpha \parallel f\parallel _{W^{1,2}\left( \Omega \right) }^{2}}{%
\mid \alpha \mid \parallel f\parallel _{W^{1,2}\left( \Omega \right) }^{2}}
\\
&=&\frac{\alpha }{\mid \alpha \mid } \\
&=&\left\{ 
\begin{array}{c}
1\text{, for }\alpha >0 \\ 
-1\text{, for }\alpha <0%
\end{array}%
\right.
\end{eqnarray*}

Thus 
\begin{equation*}
\text{\ }\alpha >0\text{ \ }\Longrightarrow \text{ \ }\theta =0\text{ \ \
and \ \ }\alpha <0\text{ }\Longrightarrow \text{ }\theta =\pi
\end{equation*}

\ \ \ \ \ \ \ \ \ \ \ \ \ \ \ \ \ \ \ \ \ \ \ \ \ \ \ \ \ \ \ \ \ \ 

\textbf{Definition 4}. The distance $\rho _{_{W^{1,2}\left( \Omega \right)
}} $ between two elements of $W^{1,2}\left( \Omega \right) $ is given by

\ \ \ 
\begin{equation*}
\rho _{_{W^{1,2}\left( \Omega \right) }}\left( f,g\right) =\text{ }\parallel
f-g\parallel _{W^{1,2}\left( \Omega \right) }
\end{equation*}

\ \ \ \ 

\textbf{Proposition 7.} For $\lambda \left( \neq 0\right) \in 
\mathbb{R}
$

\ \ 

\begin{equation*}
\rho _{_{W^{1,2}\left( \Omega \right) }}\left( f,\lambda f\right) =\mid
1-\lambda \mid \parallel f\parallel _{W^{1,2}\left( \Omega \right) }
\end{equation*}

\textbf{Proof}.

\begin{eqnarray*}
\rho _{_{W^{1,2}\left( \Omega \right) }}\left( f,\lambda f\right) &=&\left(
\dint\limits_{\Omega }\left( f-\lambda f\right) ^{2}+\left( f^{\prime
}-\lambda f^{\prime }\right) ^{2}dx\right) ^{\frac{1}{2}} \\
&=&\left( \dint\limits_{\Omega }f^{2}\left( 1-\lambda \right) ^{2}+f^{\prime
2}\left( 1-\lambda \right) ^{2}dx\right) ^{\frac{1}{2}} \\
&=&\left( \left( 1-\lambda \right) ^{2}\right) ^{\frac{1}{2}}\left(
\dint\limits_{\Omega }f^{2}+f^{\prime 2}dx\right) ^{\frac{1}{2}} \\
&=&\mid 1-\lambda \mid \parallel f\parallel _{W^{1,2}\left( \Omega \right) }
\end{eqnarray*}

\bigskip

\textbf{Corollary 1}. For $\lambda \left( \neq 0\right) \in 
\mathbb{R}
$ only when $\lambda =2$ that

\ 
\begin{eqnarray*}
\rho _{_{W^{1,2}\left( \Omega \right) }}\left( f,\lambda f\right) &=&\rho
_{_{W^{1,2}\left( \Omega \right) }}\left( f,2f\right) \\
&=&\parallel f\parallel _{W^{1,2}\left( \Omega \right) }
\end{eqnarray*}

\ 

\textbf{Corollary 2}. For $\lambda <0$ and $\lambda >2$

\begin{equation*}
\rho _{_{W^{1,2}\left( \Omega \right) }}\left( f,\lambda f\right) >\parallel
f\parallel _{W^{1,2}\left( \Omega \right) }
\end{equation*}

and for $0<\lambda <2$ we have 
\begin{equation*}
\rho _{_{W^{1,2}\left( \Omega \right) }}\left( f,\lambda f\right) <\parallel
f\parallel _{W^{1,2}\left( \Omega \right) }
\end{equation*}

\textbf{Examples 6.}

\ 

$\left( a\right) $ 
\begin{equation*}
\rho _{_{W^{1,2}\left( \Omega \right) }}\left( \cos x,\sin x\right) =\sqrt{2}
\end{equation*}

\ \ 

$\left( b\right) $ 
\begin{equation*}
\rho _{_{W^{1,2}\left( \Omega \right) }}\left( e^{x},e^{-x}\right) =\frac{%
\sqrt{e^{4}-1}}{e}
\end{equation*}

\textbf{\ \ \ \ \ \ }

\textbf{\bigskip Expansion of space: }From \textit{Proposition }$3$ above
and from the following two propositions, we see a Sobolev space is expanding
as its regularity increases.

\ \ \ \ \ \ \ 

In the next proposition we see that non zero elements are going further
apart in $W^{1,2}\left( \Omega \right) $ than they were in the Hilbert space 
$\tciLaplace ^{2}\left( \Omega \right) $

\ \ \ \ \ \ \ \ \ \ \ \ \ \ \ \ \ 

\textbf{Proposition 8}. For $f,g\in W^{1,2}\left( \Omega \right) $, 
\begin{equation*}
\rho _{_{\tciLaplace ^{2}\left( \Omega \right) }}\left( f,g\right) \leq \rho
_{_{W^{1,2}\left( \Omega \right) }}\left( f,g\right)
\end{equation*}

\ 

\textbf{Example 7}.

\ \ 

$\left( a\right) $ 
\begin{eqnarray*}
\rho _{_{W^{1,2}\left( \Omega \right) }}\left( \cos x,\sin x\right) &=&\sqrt{%
2} \\
&\geq &\sqrt{1-\sin ^{2}1} \\
&=&\rho _{_{\tciLaplace ^{2}\left( \Omega \right) }}\left( \cos x,\sin
x\right)
\end{eqnarray*}

\ \ 

$\left( b\right) $ 
\begin{eqnarray*}
\rho _{_{W^{1,2}\left( \Omega \right) }}\left( e^{x},e^{-x}\right) &=&\frac{%
\sqrt{e^{4}-1}}{e} \\
&\geq &\frac{\sqrt{e^{4}-2e^{2}-1}}{e\sqrt{2}} \\
&=&\rho _{\tciLaplace ^{2}\left( \Omega \right) }\left( e^{x},e^{-x}\right)
\end{eqnarray*}

\ \ \ \ \ \ 

\textbf{Proposition 9} (Generalizing Proposition $7$) The Sobolev space is
expanding with regularity

\ \ \ \ \ \ \ \ \ \ \ \ \ 

For $f,g\in W^{k,2}\left( \Omega \right) $ and $k\in 
\mathbb{Z}
^{+},$ we have 
\begin{equation*}
\rho _{_{W^{(k-1),2}\left( \Omega \right) }}\left( f,g\right) \leq \rho
_{_{W^{k,2}\left( \Omega \right) }}\left( f,g\right)
\end{equation*}

\textbf{Proof}. \ 
\begin{eqnarray*}
\rho _{_{W^{k,2}\left( \Omega \right) }}\left( f,g\right) &=&\left(
\dsum\limits_{j=0}^{k}\left( f^{(j)}-g^{(j)}\right) ^{2}\right) ^{\frac{1}{2}%
} \\
&\geq &\left( \dsum\limits_{j=0}^{k-1}\left( f^{(j)}-g^{(j)}\right)
^{2}\right) ^{\frac{1}{2}} \\
&=&\rho _{_{W^{(k-1),2}\left( \Omega \right) }}\left( f,g\right)
\end{eqnarray*}

\ \ \ 

where $f^{\left( j\right) }=\frac{d^{j}}{dx^{j}}\left( f\right) $ \ 

\bigskip\ \ \ \ 

\textbf{Proposition 10} (Orthogonal Decomposition).

\begin{equation*}
W^{1,2}\left( \Omega \right) =A^{2,2}\left( \Omega \right) \oplus
D^{2}\left( W_{0}^{3,2}\left( \Omega \right) \right)
\end{equation*}

where%
\begin{equation*}
A^{2,2}\left( \Omega \right) =\ker D^{2}\left( \Omega \right) \cap
W^{1,2}\left( \Omega \right)
\end{equation*}

is the Bergman space in one dimension.

\ \ 

\textbf{Proof}.

\ 

Let $f\in W^{1,2}\left( \Omega \right) $ and let 
\begin{equation*}
\eta =D_{0}^{-4}\left( D^{2}f\right)
\end{equation*}

\ 

Define 
\begin{equation*}
g=f-D^{2}\eta
\end{equation*}

Then clearly $g\in \ker D^{2}\left( \Omega \right) $ and $\eta \in
W_{0}^{3,2}\left( \Omega \right) $

\begin{equation*}
\therefore \ \ \ \ \ \ f=g\uplus D^{2}\eta
\end{equation*}

\ 

As usual there are orthogonal projections $P$ and $Q$ with

\begin{equation*}
P:W^{1,2}\left( \Omega \right) \longrightarrow A^{2,2}\left( \Omega \right)
\ \ \ \text{and \ \ }Q:W^{1,2}\left( \Omega \right) \longrightarrow
D^{2}\left( W_{0}^{3,2}\left( \Omega \right) \right)
\end{equation*}

such that 
\begin{equation*}
P+Q=I
\end{equation*}

where $I$ is the identity operator.

\ 

\textbf{Corollary 3}. 
\begin{equation*}
PQ=QP=0
\end{equation*}

and%
\begin{equation*}
P^{2}=P\text{ \ and \ }Q^{2}=Q
\end{equation*}

\textbf{\ }

\textbf{Proof}. Clearly 
\begin{equation*}
PQ=0=QP
\end{equation*}

But then the other two follow from this and the fact that 
\begin{equation*}
P+Q=I
\end{equation*}

\ 

\textbf{Example 8}. I present few but fundamntal decompositions of
elementary functions

\ 

$\left( a\right) $ For $f(x)=x,$ $P(f)=f$ and $Q(f)=0$ so that%
\begin{equation*}
f=f\uplus 0
\end{equation*}

$\left( b\right) $ For $f\left( x\right) =x^{2}$ 
\begin{equation*}
P\left( f\right) =x-\frac{1}{6}\text{ \ \ and \ \ }Q(f)=x^{2}-x+\frac{1}{6}
\end{equation*}

\ \ \ so that 
\begin{equation*}
x^{2}=\left( x-\frac{1}{6}\right) \uplus \left( x^{2}-x+\frac{1}{6}\right)
\end{equation*}

\ \ \ \ \ \ \ \ \ \ \ \ \ \ \ \ \ \ \ \ \ \ \ \ 

$\left( c\right) $ For the monomial function $f(x)=x^{n}$ 
\begin{equation*}
P(f)=\frac{6n}{n^{2}+3n+2}x-\frac{2n-2}{n^{2}+3n+2}
\end{equation*}
and 
\begin{equation*}
Q(f)=x^{n}-\frac{6n}{n^{2}+3n+2}x+\frac{2n-2}{n^{2}+3n+2}
\end{equation*}

so that \ \ 

\ \ 
\begin{equation*}
\text{ \ \ \ \ \ \ \ \ \ \ \ }x^{n}=\left( \frac{6n}{n^{2}+3n+2}x-\frac{2n-2%
}{n^{2}+3n+2}\right) \uplus \left( x^{n}-\frac{6n}{n^{2}+3n+2}x+\frac{2n-2}{%
n^{2}+3n+2}\right)
\end{equation*}

\ \ 

\ \ \ \ \ \ \ \ \ \ \ \ \ \ \ \ \ \ \ \ \ \ 

$\left( d\right) $ For $f(x)=\cos x,$

$\ \ \ \ $%
\begin{equation*}
P(f)=\left( -12+6\sin 1+12\cos 1\right) x+6-6\cos 1-2\sin 1
\end{equation*}

and

\begin{equation*}
Q(f)=\left( \cos x+(12-12\cos 1-6\sin 1\right) x-6+6\cos 1+2\sin 1
\end{equation*}

so that%
\begin{equation*}
f=P(f)\uplus Q(f)
\end{equation*}

\ \ \ \ \ \ \ 

$\left( f\right) $ The last example I consider is $f(x)=e^{x}$

\begin{equation*}
e^{x}=\underset{P\left( f\right) }{\underbrace{\left( -6ex+4e\right) }}%
\uplus \underset{Q(f)}{\underbrace{\left( e^{x}+6ex-4e\right) }}
\end{equation*}

Indeed, let 
\begin{equation*}
\eta =D_{0}^{-4}\left( D^{2}f\right) =D_{0}^{-4}\left( e^{x}\right)
\end{equation*}

up on solving the differential equation with vanishing boundary conditions
inversely, we have

\ \ \ 
\begin{equation*}
\eta \left( x\right) =e^{x}+ex^{3}-2ex^{2}.
\end{equation*}

Consider 
\begin{equation*}
g:=f-D^{2}\eta =-6ex+4e
\end{equation*}%
so that 
\begin{equation*}
f=g\uplus D^{2}\eta
\end{equation*}%
That is 
\begin{equation*}
e^{x}=\left( -6ex+4e\right) \uplus \left( e^{x}+6ex-4e\right)
\end{equation*}

Up on calculations, we see 
\begin{equation*}
P\left( -6ex+4e\right) =-6ex+4e\text{ \ and \ \ }Q\left( -6ex+4e\right) =0
\end{equation*}
justifying the fact that 
\begin{equation*}
P^{2}=P\text{\ \ \ and \ \ \ }Q\circ P=0
\end{equation*}

\ \ \ \ 

\ 

\textbf{Proposition 11}. For $f\in W^{1,2}\left( \Omega \right) $, we have 
\begin{equation*}
\left\langle Pf,Qf\right\rangle _{\tciLaplace ^{2}\left( \Omega \right)
}=-\left\langle \left( Pf\right) ^{\prime },\left( Qf\right) ^{\prime
}\right\rangle _{\tciLaplace ^{2}\left( \Omega \right) }
\end{equation*}

\ 

\textbf{Proof}. This follows from the fact that $Pf$ and $Qf$ are orthogonal
in the $W^{1,2}\left( \Omega \right) $.

\ \ \ \ 

\textbf{Example 9.}

\ 

$\left( a\right) $ 
\begin{equation*}
\left\langle x-\frac{1}{6},x^{2}-x+\frac{1}{6}\right\rangle _{\tciLaplace
^{2}\left( \Omega \right) }=-\left\langle 1,2x-1\right\rangle _{\tciLaplace
^{2}\left( \Omega \right) }
\end{equation*}

$\left( b\right) $ For $\alpha \beta =-1$, 
\begin{equation*}
\left\langle e^{\alpha x},e^{\beta x}\right\rangle _{\tciLaplace ^{2}\left(
\Omega \right) }=-\left\langle \alpha e^{\alpha x},\beta e^{\beta
x}\right\rangle _{\tciLaplace ^{2}\left( \Omega \right) }
\end{equation*}

\ \ \ \ \ \ \ \ \ \ \ \ \ \ \ \ \ \ \ \ \ \ 

\textbf{The space }$W^{-1,2}\left( \Omega \right) $\ \ 

\ \ \ \ 

The space $W^{-1,2}\left( \Omega \right) $\ is the conjugate space of the
Sobolev space $W^{1,2}\left( \Omega \right) $, that is

\ \ 
\begin{equation*}
W^{-1,2}\left( \Omega \right) =\left( W^{1,2}\left( \Omega \right) \right)
^{\ast }
\end{equation*}

where 
\begin{equation*}
\left( W^{1,2}\left( \Omega \right) \right) ^{\ast }:=\left\{ \tau
:W^{1,2}\left( \Omega \right) \longrightarrow 
\mathbb{R}
,\tau \text{ is a bounded linear functional}\right\}
\end{equation*}

The spaces $W_{0}^{1,2}\left( \Omega \right) $ and $W_{0}^{1,2}\left( \Omega
\right) ^{\perp }$

\bigskip

We have 
\begin{equation*}
W_{0}^{1,2}\left( \Omega \right) ^{\perp }=\left\{ f\in W^{1,2}\left( \Omega
\right) :\langle f,g\rangle =0,\forall g\in W_{0}^{1,2}\left( \Omega \right)
\right\}
\end{equation*}

\ \ 

and from linear algebra of vector spaces we have a direct sum

\begin{equation*}
W^{1,2}\left( \Omega \right) =W_{0}^{1,2}\left( \Omega \right) \oplus
W_{0}^{1,2}\left( \Omega \right) ^{\perp }
\end{equation*}

\bigskip and therefore we have an interesting relation 
\begin{equation*}
W^{1,2}\left( \Omega \right) \ominus W_{0}^{1,2}\left( \Omega \right)
=W_{0}^{1,2}\left( \Omega \right) ^{\perp }
\end{equation*}%
\ \ \ 
\begin{equation*}
W^{1,2}\left( \Omega \right) \ominus W_{0}^{1,2}\left( \Omega \right)
\end{equation*}

\ 

We next see that in fact the above direct difference or simply $%
W_{0}^{1,2}\left( \Omega \right) ^{\perp }$ is a two dimensional subspace of 
$W^{1,2}\left( \Omega \right) $

\ 

Indeed

\begin{equation*}
W_{0}^{1,2}\left( \Omega \right) ^{\perp }=\left\{ f\in W^{1,2}\left( \Omega
\right) :\langle f,g\rangle =0,\forall g\in W_{0}^{1,2}\left( \Omega \right)
\right\}
\end{equation*}

where 
\begin{equation*}
\langle f,g\rangle =\int_{\Omega }\left( fg+f^{\prime }g^{\prime }\right)
dx=0
\end{equation*}

\ 

Applying integration by parts with no boundary integrals as $g\in
W_{0}^{1,2}\left( \Omega \right) ,$ we have

\ \ 

\begin{eqnarray*}
\int_{\Omega }\left( fg-f^{\prime ^{\prime }}g\right) dx &=&\int_{\Omega
}\left( f-f^{\prime ^{\prime }}\right) gdx \\
&=&0,\text{ \ \ \ }\forall g\in W_{0}^{1,2}\left( \Omega \right) \\
\text{ \ \ \ } &\Longrightarrow &\text{\ \ \ }f-f^{\prime \prime }=0
\end{eqnarray*}

\ \ \ 

Solving the second order ordinary differential equation $f^{\prime \prime
}-f=0$ we have

\begin{equation*}
f_{c}=\alpha e^{x}+\beta e^{-x}\in W_{0}^{1,2}\left( \Omega \right) ^{\perp }
\end{equation*}
\ 

for $\alpha ,\beta $ arbitrary real constants.

\ \ 

Therefore we have 
\begin{equation*}
\left( W_{0}^{1,2}\left( \Omega \right) \right) ^{\perp }=\text{span}\langle
e^{x},e^{-x}\rangle
\end{equation*}

From the fact that

\ 
\begin{equation*}
W^{1,2}\left( \Omega \right) =W_{0}^{1,2}\left( \Omega \right) \oplus \left(
W_{0}^{1,2}\left( \Omega \right) \right) ^{\perp }
\end{equation*}

we have

\ \ 
\begin{equation*}
W^{1,2}\left( \Omega \right) \ominus W_{0}^{1,2}\left( \Omega \right)
=\left( W_{0}^{1,2}\left( \Omega \right) \right) ^{\perp }=\text{span}%
\langle e^{x},e^{-x}\rangle
\end{equation*}

\ \ 

Hence $W^{1,2}\left( \Omega \right) \ominus W_{0}^{1,2}\left( \Omega \right) 
$ is a skinny two dimensional subspace. This is interesting by it self,
showing the fact that when we remove all elements that vanish on the
boundary, the space what is left is a two dimensional subspace.

\ 

From the direct sum 
\begin{equation*}
W^{1,2}\left( \Omega \right) =W_{0}^{1,2}\left( \Omega \right) \oplus \left(
W_{0}^{1,2}\left( \Omega \right) \right) ^{\perp }
\end{equation*}

we have orthogonal projections 
\begin{equation*}
\overset{\sim }{P}:W^{1,2}\left( \Omega \right) \longrightarrow
W_{0}^{1,2}\left( \Omega \right) \text{ \ \ and \ }\overset{\sim }{Q}%
:W^{1,2}\left( \Omega \right) =W_{0}^{1,2}\left( \Omega \right) ^{\perp }
\end{equation*}

such that $\forall f\in W^{1,2}\left( \Omega \right) ,$ 
\begin{equation*}
\overset{\sim }{P}(f)=f-\left( \alpha e^{x}+\beta e^{-x}\right) \ ,\overset{%
\sim }{Q}(f)=\alpha e^{x}+\beta e^{-x}
\end{equation*}

where 
\begin{equation*}
\alpha =\frac{f(1)e-f(0)}{e^{2}-1}\text{ \ and \ }\beta =\frac{%
f(0)e^{2}-f(1)e}{e^{2}-1}
\end{equation*}

\bigskip

\textbf{Proposition 12. }$\overset{\sim }{P}(f)_{\mid \partial \Omega }=0$
and $\overset{\sim }{Q}f_{\mid \partial \Omega }=f_{\mid \partial \Omega }$

\ 

\textbf{Proof}. \ Clearly $\overset{\sim }{P}(f)\left( 0\right) =0$ and $%
\overset{\sim }{P}(f)\left( 1\right) =0$

\begin{eqnarray*}
\ \ \ \overset{\sim }{Q}(f)(0) &=&\alpha +\beta \\
&=&\frac{f(1)e-f(0)}{e^{2}-1}+\frac{f(0)e^{2}-f(1)e}{e^{2}-1} \\
&=&\frac{f(0)(e^{2}-1)}{e^{2}-1} \\
&=&f(0)
\end{eqnarray*}

\ \ \ \ \ \ \ \ \ \ \ \ \ \ \ 

\begin{eqnarray*}
\ \ \overset{\sim }{Q}(f)(1) &=&\alpha e+\beta e^{-1} \\
&=&\frac{f(1)e^{2}-f(0)e}{e^{2}-1}+\frac{f(0)e-f(1)}{e^{2}-1} \\
&=&\frac{f(1)\left( e^{2}-1\right) }{e^{2}-1} \\
&=&f(1)
\end{eqnarray*}

\ \ 

Now from 
\begin{equation*}
W^{1,2}\left( \Omega \right) =W_{0}^{1,2}\left( \Omega \right) \oplus \left(
W_{0}^{1,2}\left( \Omega \right) \right) ^{\perp }
\end{equation*}
we have 
\begin{equation*}
W^{1,2}\left( \Omega \right) ^{\ast }=W_{0}^{1,2}\left( \Omega \right)
^{\ast }\oplus \left( W_{0}^{1,2}\left( \Omega \right) \right) ^{\perp \ast }
\end{equation*}

That is 
\begin{equation*}
W^{-1,2}\left( \Omega \right) =W_{0}^{-1,2}\left( \Omega \right) ^{\ast
}\oplus \left( W_{0}^{1,2}\left( \Omega \right) \right) ^{\perp \ast }
\end{equation*}

Looking 
\begin{equation*}
W_{0}^{-1,2}\left( \Omega \right) ^{\ast }=\left\{ \zeta :W_{0}^{1,2}\left(
\Omega \right) \longrightarrow 
\mathbb{R}
\text{, bounded linear functional}\right\}
\end{equation*}

By Riesz representation theorem 
\begin{equation*}
\exists f\in W_{0}^{1,2}\left( \Omega \right) :\gamma \left( g\right)
=\left\langle \gamma ,g\right\rangle =\dint\limits_{\Omega }gf+g^{\prime
}f^{\prime }dx,\forall g\in W_{0}^{1,2}\left( \Omega \right)
\end{equation*}

Then we have the representation of $\gamma $ to be 
\begin{equation*}
\gamma =f-f^{\prime \prime }
\end{equation*}%
where $f,$ $f^{\prime }\in \tciLaplace ^{2}\left( \Omega \right) $.

\ \ 

If $\gamma \neq 0$, then the function $f$ that is used to represent $\gamma $
is a solution of the inhomogeneous differential equation

\ \ \ \ \ \ \ \ \ \ \ \ \ \ \ \ \ \ \ \ 
\begin{equation*}
\gamma =f-f^{\prime \prime }
\end{equation*}

\ \ 

The function that used to represent the zero linear functional $\gamma =0$
is a solution of the homogeneous ordinary differential equation given by

\bigskip 
\begin{equation*}
f^{\prime \prime }-f=0
\end{equation*}%
.

Up on solving the latter equation, we get the function that represents the
zero linear functional to be

\ 
\begin{equation*}
f_{c}=\alpha e^{x}+\beta e^{-x}
\end{equation*}%
with $\alpha ,\beta $ some real constants.

\ \ 

\textbf{Proposition 13 (}Representation of a linear functional\textbf{). }

\ 

Let $\zeta :W_{0}^{1,2}\left( \Omega \right) \longrightarrow 
\mathbb{R}
$ be a bounded linear functional, then by Riesz representation theorem, $%
\exists f\in W_{0}^{1,2}\left( \Omega \right) $ such that

\ \ $\ \ $%
\begin{equation*}
\zeta =f_{0}-\frac{d}{dx}f_{1}\ \ \ \text{with \ }\ f_{0}=f,\text{ \ }f_{1}=%
\frac{d}{dx}f=f^{\prime }
\end{equation*}%
where $f,$ $f_{1}\in \tciLaplace ^{2}\left( \Omega \right) $

\ 

\textbf{\ }

\textbf{Proof}. For $\zeta \in W_{0}^{1,2}\left( \Omega \right) ^{\ast }$ , $%
\exists f\in W_{0}^{1,2}\left( \Omega \right) :$%
\begin{eqnarray*}
\zeta \left( g\right) &=&\left\langle \zeta ,g\right\rangle \\
&=&\dint\limits_{\Omega }\left( fg+f^{\prime }g^{\prime }\right) dx,\forall
g\in W_{0}^{1,2}\left( \Omega \right)
\end{eqnarray*}

But%
\begin{eqnarray*}
\dint\limits_{\Omega }\left( fg+f^{\prime }g^{\prime }\right) dx
&=&\dint\limits_{\Omega }\left( fg-f^{^{\prime }\prime }g\right) dx \\
&=&\dint\limits_{\Omega }\left( f-f^{\prime \prime }\right) gdx \\
&=&\left\langle f-f^{\prime \prime },g\right\rangle
\end{eqnarray*}

\ \ 

\begin{equation*}
\therefore \ \ \ \zeta =f_{0}-\frac{d}{dx}f_{1}\text{ \ \ with \ \ }f_{0}=f,%
\text{ \ }f_{1}=\frac{d}{dx}f=f^{\prime }
\end{equation*}

where $f_{0,}$ $f_{1}\in \tciLaplace ^{2}\left( \Omega \right) $.

\ \ \ \ 

\textbf{Definition 5} (Projection).

\ \ \ \ 

Let $f,g:\Omega \longrightarrow 
\mathbb{R}
$ with $\parallel g\parallel _{W^{1,2}\left( \Omega \right) }\neq 0$. Then
we define the projection of $f$ over $g$ denoted

\ 
\begin{equation*}
\Pr \text{oj}_{g}\left( f\right) _{W^{1,2}\left( \Omega \right) }:=\frac{%
\left\langle f,g\right\rangle _{W^{1,2}\left( \Omega \right) }}{\parallel
g\parallel _{W^{1,2}\left( \Omega \right) }^{2}}g
\end{equation*}

\ \ \ 

\textbf{Proposition 14}. For $f,g,h:\Omega \longrightarrow 
\mathbb{R}
$ with $\parallel h\parallel _{W^{1,2}\left( \Omega \right) }\neq 0$ and $%
\alpha \in 
\mathbb{R}
$, we have

\ \ 

$\left( a\right) \Pr $oj$_{h}\left( f+g\right) _{W^{1,2}\left( \Omega
\right) }=\Pr $oj$_{h}\left( f\right) _{W^{1,2}\left( \Omega \right) }+\Pr $%
oj$_{h}\left( g\right) _{W^{1,2}\left( \Omega \right) }$

\ \ \ \ 

$\left( b\right) $ $\Pr $oj$_{h}\left( \alpha f\right) _{W^{1,2}\left(
\Omega \right) }=\alpha \Pr $oj$_{h}\left( f\right) _{W^{1,2}\left( \Omega
\right) }$

\ \ \ \ \ \ \ \ \ \ \ \ \ 

\textbf{Corollary 4.}

\ \ \ \ 

$\left( a\right) $ If the two functions are orthogonal then $\Pr $oj$%
_{g}\left( f\right) _{W^{1,2}\left( \Omega \right) }=0$

\ 

$\left( b\right) $ $\Pr $oj$_{f}(f)_{W^{1,2}\left( \Omega \right) }=f$

\ 

$\left( d\right) $ $\Pr $oj$\left( _{\beta g}\right) \left( f\right)
_{W^{1,2}\left( \Omega \right) }=\Pr $oj$_{g}\left( f\right) _{W^{1,2}\left(
\Omega \right) }$

\ 

$\left( e\right) $ $\Pr $oj$_{\left( \beta g\right) }\left( \alpha f\right)
_{W^{1,2}\left( \Omega \right) }=\alpha \Pr $oj$_{g}\left( f\right)
_{W^{1,2}\left( \Omega \right) }$

\ \ \ \ 

\textbf{Example 10.}

\ \ 

$\left( a\right) $ $\ \Pr $oj$_{\left( x^{2}\right) }e_{W^{1,2}\left( \Omega
\right) }^{x}=\frac{15e}{23}x^{2}$

\ \ \ 

$\left( b\right) $ $\ \Pr $oj$\left( _{e^{-x}}\right) e_{W^{1,2}\left(
\Omega \right) }^{x}=0$

\ \ 

$\left( c\right) $ \ $\Pr $oj$\left( _{\cos x}\right) \left( \sin x\right)
_{W^{1,2}\left( \Omega \right) }=0$

\ \ \ \ 

$\left( d\right) $ \ $\Pr $oj$\left( _{e^{\beta x}}\right) e_{W^{1,2}\left(
\Omega \right) }^{\alpha x}=\gamma e^{\beta x}$ for $\alpha ,\beta \neq 0$
and $\gamma =\frac{\left( 2\alpha +2\beta \right) e^{\alpha +\beta }-2\alpha
\beta ^{2}-2\beta }{\left( \alpha +\beta \right) \left( \beta ^{2}+1\right)
\left( e^{2\beta }-1\right) }$

\ \ \ \ 

\ 

\ \ 

\textbf{References}

\ 

\ 

$\left[ 1\right] $ \ Dejenie A. Lakew, New Proofs on Properties of an
Orthogonal Decomposition

\ \ \ \ 

\ \ \ \ \ \ of a Hilbert Space, $arXiv:1510.07944v1$

\ 

$\left[ 2\right] $ \ Dejenie A. Lakew, On Orthogonal Decomposition of $%
\tciLaplace 
{{}^2}%
(%
\Omega
)$, Int. J. Math.

\ \ \ \ \ Comp. Sc., $10(2015),$No$.1,$ $27-37$

\end{document}